\documentclass[letterpaper,12pt,reqno]{amsart}

\usepackage[margin=1.5in]{geometry}
\usepackage{hyperref}

\hypersetup{
     colorlinks   = true,
     citecolor    = black,
     linkcolor    = blue
} 
 
\usepackage{amsmath}
\usepackage{amssymb}
\usepackage{marginnote}

\usepackage{enumitem}
\usepackage{color}
\usepackage{url}
\usepackage[usenames,dvipsnames]{xcolor}
\usepackage{graphicx}
\usepackage{verbatim}
\usepackage{tikz}

\makeatletter 
\def\@cite#1#2{{\m@th\upshape\bfseries%
[{#1\if@tempswa{\m@th\upshape\mdseries, #2}\fi}]}}
\makeatother 

\theoremstyle{plain}
\newtheorem{thm}{Theorem}[section]

\newtheorem{cor}[thm]{Corollary}
\newtheorem{prop}[thm]{Proposition}

\theoremstyle{definition}

\newtheorem{quest}[thm]{Question}

\theoremstyle{}

\numberwithin{equation}{subsection}

\newcommand{\nc}{\newcommand}
\newcommand{\rnc}{\renewcommand}

\nc\bA{\mathbb{A}}
\nc\bB{\mathbb{B}}
\nc\bC{\mathbb{C}}
\nc\bD{\mathbb{D}}
\nc\bE{\mathbb{E}}
\nc\bF{\mathbb{F}}
\nc\bG{\mathbb{G}}
\nc\bH{\mathbb{H}}
\nc\bI{\mathbb{I}}
\nc{\bJ}{\mathbb{J}} 
\nc\bK{\mathbb{K}}
\nc\bL{\mathbb{L}}
\nc\bM{\mathbb{M}}
\nc\bN{\mathbb{N}}
\nc\bO{\mathbb{O}}
\nc\bP{\mathbb{P}}
\nc\bQ{\mathbb{Q}}
\nc\bR{\mathbb{R}}
\nc\bS{\mathbb{S}}
\nc\bT{\mathbb{T}}
\nc\bU{\mathbb{U}}
\nc\bV{\mathbb{V}}
\nc\bW{\mathbb{W}}
\nc\bY{\mathbb{Y}}
\nc\bX{\mathbb{X}}
\nc\bZ{\mathbb{Z}}
\nc\wt{\widetilde}

\nc\cA{\mathcal{A}}
\nc\cB{\mathcal{B}}
\nc\cC{\mathcal{C}}
\nc\cD{\mathcal{D}}
\nc\cE{\mathcal{E}}
\nc\cF{\mathcal{F}}
\nc\cG{\mathcal{G}}
\nc\cH{\mathcal{H}}
\nc\cI{\mathcal{I}}
\nc{\cJ}{\mathcal{J}} 
\nc\cK{\mathcal{K}}
\nc\cL{\mathcal{L}}
\nc\cM{\mathcal{M}}
\nc\cN{\mathcal{N}}
\nc\cO{\mathcal{O}}
\nc\cP{\mathcal{P}}
\nc\cQ{\mathcal{Q}}
\nc\cR{\mathcal{R}}
\nc\cS{\mathcal{S}}
\nc\cT{\mathcal{T}}
\nc\cU{\mathcal{U}}
\nc\cV{\mathcal{V}}
\nc\cW{\mathcal{W}}
\nc\cY{\mathcal{Y}}
\nc\cX{\mathcal{X}}
\nc\cZ{\mathcal{Z}}

\newcommand{\defeq}{:=}

\nc{\Tm}{Teichm\"uller}
\nc{\Mod}{\text{Mod}}

\newcommand{\SL}{\mathrm{SL}}

\newcommand{\SO}{\mathrm{SO}}

\newcommand{\bfV}{\textbf{V}}
\newcommand{\bfv}{\textbf{v}}

\rnc{\Re}{\text{Re}}
\rnc{\Im}{\text{Im}}
\nc{\var}{\text{var}}
\nc{\cov}{\text{covar}}
\nc{\deq}{\stackrel{d}{=}}

\def\strutdepth{\dp\strutbox}
\def \ss{\strut\vadjust{\kern-\strutdepth \sss}}
\def \sss{\vtop to \strutdepth{
\baselineskip\strutdepth\vss\llap{$\diamondsuit\;\;$}\null}}

\def\strutdepth{\dp\strutbox}
\def \sst{\strut\vadjust{\kern-\strutdepth \ssss}}
\def \ssss{\vtop to \strutdepth{
\baselineskip\strutdepth\vss\llap{$\spadesuit\;\;$}\null}}

\def\strutdepth{\dp\strutbox}
\def \ssh{\strut\vadjust{\kern-\strutdepth \sssh}}
\def \sssh{\vtop to \strutdepth{
\baselineskip\strutdepth\vss\llap{$\heartsuit\;\;$}\null}}

\nc{\dmo}{\DeclareMathOperator}
\dmo{\rank}{rank}
\dmo{\End}{End}
\dmo{\Hom}{Hom}
\dmo{\Jac}{Jac}
\dmo{\Id}{Id}
\dmo{\Ann}{Ann}
\dmo{\Area}{Area}
\dmo{\CP}{\bC P^1}
\dmo{\Aut}{Aut}

\title{Oscillations and the Kontsevich-Zorich cocycle}

\author[H. Al-Saqban]{Hamid Al-Saqban}
\address{Department of Mathematics, University of Maryland, College Park, MD}
\email{hqs@math.umd.edu}

\begin{document}
\maketitle
\begin{abstract}
We present a mechanism for producing oscillations along the lift of the Teichm\"uller geodesic flow to the (real) Hodge bundle, as the
basepoint surface is deformed by a unipotent element of $\SL_2(\bR)$. Invoking Chen-M\"oller \cite{CM13}, we apply our methods to all but finitely many strata in genus $4$, those exhibiting a varying Lyapunov-exponents phenomenon. 
\end{abstract}

\section{Introduction}

Let $\pi: \textbf{H} \to X$ be the absolute (real) Hodge bundle over an $\SL_2(\bR)$ orbit closure $X$, whose $2g$-dimensional fiber over each point in $X$ is $H^1(S,\bR)$. Let $\nu$ be an ergodic $\SL_2(\bR)$-invariant probability measure on $X$. For $g\in\SL_2(\bR)$, the Kontsevich-Zorich cocycle $g_\ast$ is the lift of the action of $g$ to $\textbf{H}$, obtained by parallel transport with respect to the Gauss-Manin connection. Moreover, $g_\ast$ acts symplectically since it preserves the intersection form on $H^1(S,\bR)$. 

Let $h_t$ be the \Tm~horocycle flow and let $\text{T}h_t \defeq (h_t)_\ast$ be its lift to the projectivized Hodge bundle $\mathbb{P}(\mathbf{H})$. Suppose that one is interested in studying the probability measures on $\mathbb{P}(\mathbf{H})$ that are invariant under $\text{T}h_t$. Motivated by the work of Bainbridge-Smillie-Weiss \cite{BSW16}, Forni posed the following question:

\begin{quest}[G. Forni]
Let $\hat{\mu}$ be a $\text{T}h_t $-invariant probability measure supported on the projectivized bundle $\mathbb{P}(\mathbf{H})$ and which is not supported on the zero-section. Is it true that the push-forward measure $\mu$ under the projection map $\pi$ to the moduli space must be supported on an orbit closure with completely degenerate Kontsevich-Zorich exponents?
\end{quest}

In the symplectic orthogonal of the tautological subbundle, it is known that orbit closures with completely degenerate Kontsevich-Zorich exponents exist in the strata $\cH(1,1,1,1)$ and $\cH^{\text{even}}(2,2,2)$ \cite{For06, FM08,FMZ11}, and are referred to in the literature as \emph{Eierlegende Wollmilchsau} and \emph{Ornithorynque}. That these are the only orbit closures with completely degenerate Kontsevich-Zorich exponents follows from the works \cite{Mol05, EKZ14, Aul16, AN19}. It is also known that the cocycle acts by isometries in these two orbit closures, and so by a Krylov-Bogoliubov construction, one can construct $\text{T}h_t $-invariant probability measures that are non-trivial (that is, not supported on the zero-section). One can then ask whether these are the only possible examples that arise. In all other situations, we expect that oscillations of the norm of the cohomology classes along large circles, when they exist, prevent the existence of non-trivial $\text{T}h_t $-invariant probability measures in general.

In this paper, we give partial evidence (in the affirmative) towards Forni's question, and present a mechanism for producing oscillations along the lift of the Teichm\"uller geodesic flow to the (real) Hodge bundle, as the basepoint surface is rotated. We apply our methods to all but finitely many strata in genus 4, those exhibiting a varying Lyapunov-exponents phenomenon, as is shown in Chen-M\"oller \cite{CM13}. 

More precisely, fix a norm $\|\cdot\|_{\pi(\cdot)}$ on $\textbf{H}$. Define $\sigma : \SL_2(\bR) \times \textbf{H} \to \bR$ by \[\sigma(g,\textbf{v}) = \frac{\|g_{\ast} \textbf{v}\|_{g\pi(\textbf{v})}}{\|\textbf{v}\|_{\pi(\bfv)}}\]

Let $\textbf{V}$\label{V} be a $\nu$-strongly irreducible $\SL_2(\bR)$-invariant subbundle in the symplectic orthogonal of the tautological subbundle, which is defined and is continuous on $X$ \cite{Fil06}.

For each $\omega \in X$ with full orbit closure, and $\textbf{v}_\omega\neq 0$ in $\text{\bf{V}}_\omega$, and for a.e. $r_\theta \in \SO_2(\bR)$, it is a consequence of a theorem of Chaika-Eskin \cite{CE13} that \[\lim_{t\to\infty}\frac{1}{t}\log\sigma(g_t r_\theta, \textbf{v}_\omega)  = \lambda\] where $\lambda = \lambda(\nu)$ is the top Lyapunov exponent of the restriction of the Kontsevich-Zorich cocycle to $\textbf{V}$. See \cite{AAE+17} for a  refinement of this result to a positive Hausdorff codimension set of angles. 

We say that an orbit closure $X$ has a \textbf{varying Lyapunov phenomenon} if there exists an affine invariant submanifold of $X$ that supports an ergodic $\SL_2(\bR)$-invariant measure $\nu'$ and such that $\lambda(\nu') \neq \lambda(\nu)$. For such an  $X$, we show that the cocycle cannot be normalized by a function that is independent of $\theta$, for a full measure set of angles. 

More precisely, letting $h_s=\begin{pmatrix} 1&s\\0&1\end{pmatrix}, \, \bar{h}_s=\begin{pmatrix} 1&0\\s&1 \end{pmatrix}, \, g_t=\begin{pmatrix}e^{t}&0\\0&e^{-t}\end{pmatrix} $ and $r_\theta=\begin{pmatrix} \cos(\theta)&-\sin(\theta)\\ \sin(\theta)&\cos(\theta) \end{pmatrix}$, we show 

\begin{thm}\label{maintheta}
For any $\omega$ in $X$ so that $\overline{\SL_2(\mathbb{R}) \cdot \omega}=X$ with $\lambda = \lambda(\nu)>0$, and such that $X$ has a varying Lyapunov phenomenon, and for any function $f(t)$, the set \[\left\{r_\theta \in \SO_2(\bR) ~:~ \lim_{t\to\infty} \frac{\sigma(g_t r_\theta, \text{\bf{v}}_\omega) }{f(t)} \text{ converges to a non-zero number}\right\}\] has zero measure with respect to the Haar measure on $\SO_2(\bR)$. 
\end{thm}

We show that Theorem \ref{maintheta} follows from its horocyclic counterpart

\begin{thm}\label{mainhoro}
Under the hypothesis in Theorem \ref{maintheta}, the set \[\left\{s \in [-1,1] ~:~ \lim_{t\to\infty} \frac{\sigma(g_t h_s, \text{\bf{v}}_\omega) }{f(t)} \text{ converges to a non-zero number}\right\}\] has zero measure with respect to the Lebesgue measure on $[-1,1]$. 
\end{thm}

\subsection{Sketch of the Proof} The idea of the proof is to construct two regimes, one in which the Kontsevich-Zorich cocycle grows as expected, and another where the behavior is atypical. 

The expected behavior is an input that comes from the work of Chaika-Eskin \cite{CE13}, where it is shown that one can construct a large open set in which the cocycle grows as expected in finite time (i.e. exponentially with rate $\lambda$ up to some additive error), and whose set of bad futures is small.

The atypical behavior is extracted from carefully defined tubular neighborhoods of a submanifold of $X$ whose (second) largest Lyapunov exponent is varying. 
To conclude the argument, we first use quantitive recurrence results from Eskin-Mirzakhani-Mohammadi \cite{EMM15} to control returns to the typical and atypical neighborhoods, then apply a Lebesgue density argument to deduce that the limit cannot exist for a full measure set of directions. 

\subsection{Applications}
For $\omega\in X$, let $\textbf{v}_\omega$ in $\textbf{H}_\omega$ be any isotropic subspace of dimension $g-1$ in the symplectic orthogonal of the tautological subspace, $\text{span}\{[\text{Re} ~\omega], [\text{Im} ~ \omega]\}$, of the Hodge bundle $\textbf{H}$. For each $\omega \in X$ with full orbit closure, and for a.e. $r_\theta \in \SO_2(\bR)$, it is also a consequence of the same theorem of Chaika-Eskin \cite{CE13} referred to above that \[\lim_{t\to\infty} \frac{1}{t}\log\sigma(g_t, \textbf{v}_\omega) = \sum_{i=2}^g \lambda_i\] where, together with $\lambda_1 = 1$, $\lambda_i$ are the top $g$ Lyapunov exponents of the Kontsevich-Zorich cocycle. The top $g$ exponents determine the entire Lyapunov spectrum by symplecticity. The work of Chen-M\"oller \cite{CM13} demonstrates, in part, the following result:

\begin{thm}\cite[Theorem 1.2]{CM13}\label{nonvarying}
For all but finitely many strata in genus 4, the sum of Lyapunov exponents is varying.
\end{thm}

Thus the conclusions of Theorem \ref{maintheta} and \ref{mainhoro} hold for any surface with full orbit closure in all but finitely many strata in genus 4.

\subsection{Related results} In \cite{DFV17}, Dolgopyat-Fayad-Vinogradov prove a central limit theorem for the time integral of sufficiently regular zero-average observables of the pushforward of a small horocyclic arc by the geodesic flow, as the basepoint varies generically with respect to an ergodic $P$-invariant measure, where $P$ is the upper triangular subgroup of $\SL_2(\bR)$ (their results are in fact much more general, cf. \cite[Theorem 7.1]{DFV17}, but we present their theorem in the $\SL_2(\bR)$ setting for simplicity). While their results are not immediately applicable to the Kontsevich-Zorich cocycle, it shows that Theorem \ref{mainhoro}, which is, in our specialized setting, a qualitative analogue of their results, can likely be strengthened quantitatively, if one allows the surface to vary generically in an orbit closure (so this would be a tradeoff, since our result is true for any fixed basepoint, thanks to the work of Eskin-Mirzakhani-Mohammadi). That a central limit theorem holds for the Kontsevich-Zorich cocycle as a basepoint varies generically is the subject of \cite{A19}, and the approach in that paper uses completely different tools (the Brownian motion).

Recently, and after the completion of this paper, the oscillation mechanism presented here has been greatly developed and refined by Chaika, Khalil, and Smillie in their work on limit measures of \Tm~horocyclic arcs (cf. \cite[Theorem 4.1]{CKS21}).
 
\section*{Acknowledgements} 
I am grateful to Jon Chaika for hosting me in Salt Lake City where this work was begun in July 2016, for explaining the tools used in this paper, and for numerous helpful conversations, without which this paper would not exist. I am grateful to my advisor Giovanni Forni for asking the question that led to this project, and for his guidance and support. I gratefully acknowledge support from the NSF grants DMS 1107452, 1107263, 1107367 ``RNMS: Geometric Structures and Representation Varieties" (the GEAR Network) and DMS 1600687. 

  \section{Preliminaries}
  \subsection{Translation surfaces}
Let $S$ be a Riemann surface of genus $g\geq 2$, and $\omega$ a holomorphic $1$-form on $S$. The pair $(S,\omega)$ is said to be a translation surface, since $\omega$ gives a (degenerate) flat metric on $S$, and $\omega$ is invariant under translations when it is written in local coordinates. The zero set $\Sigma$ of $\omega$ characterizes the singularity set of the conical metric. The area of a translation surface is given by $\int_S \omega \wedge \overline{\omega}$. We refer to the pair $(S,\omega)$ as just $\omega$. 

\subsection{Moduli Space}
Let $\cT\cH_g$ be the \Tm~space of unit-area translation surfaces of genus $g \geq 2$, and let $\cH_g = \cT\cH_g / \Mod_g$ be the corresponding moduli space, where $\Mod_g$ is the mapping class group. The space $\cH_g$ is partitioned into strata $\cH(\kappa) = \cH(\kappa_1,\ldots,\kappa_n)$, which consist of unit-area translation surfaces whose singularities have cone angle $2\pi\kappa_i$, and $\sum \kappa_i = 2g-2$. One can also define local period coordinates in a stratum, where all changes of coordinates are given by affine maps. 

\subsection{$\SL_2(\bR)$ action}
There is a natural action of $\SL_2(\bR)$ on translation surfaces and on their moduli. It is shown in \cite{EM18, EMM15} that for any $\omega \in \cH(\kappa)$, the closure $X$ of $\SL_2(\bR) \cdot \omega$ is an affine invariant submanifold, and supports an ergodic $\SL_2(\bR)$-invariant probability measure $\nu$.

\subsection{Hodge inner product} Given two holomorphic $1$-forms $\omega_1,\omega_2$ in $\Omega(S)$, where $\Omega(S)$ is the vector space of holomorphic $1$-forms on $S$, the Hodge inner product is defined to be \[\langle\omega_1,\omega_2 \rangle := \frac{i}{2}\int_S \omega_1 \wedge \overline{\omega_2}\]
Moreover, the Hodge representation theorem implies that for any given cohomology class $c \in H^1(S,\bR)$, there is a unique holomorphic $1$-form $h(c) \in \Omega(S)$, such that $c = [\Re~ h(c)]$. We define the Hodge inner product for two real cohomology classes $c_1,c_2 \in H^1(S,\bR)$ as \[(c_1,c_2)_{\omega} := \langle h(c_1),h(c_2) \rangle_{\omega}\]

For $c$ in the symplectic orthogonal of $[\omega]$, it also follows from the work of Forni \cite{For02} (see also \cite[Corollary 2.1]{FMZ12} and \cite[Corollary 30]{FM13}) that  \begin{align}\label{lipschitz} \left|\frac{d}{dt} \log \|c\|_{g_t \omega}\right| < 1\end{align} 

\subsection{Expected behavior}\label{expected}
For each of the ambient manifold $X$ and a submanifold $X'$, we will need to construct open sets where the cocycle grows as expected up to some additive error, and whose set of bad futures is small. This is implemented in \cite{CE13} for random walks, and is adapted to the deterministic case in \cite[Corollary 7.6]{AAE+17}. The main point is that random walks track geodesics up to an error that is sublinear in hyperbolic distance, and the rest of the adaptation follows by standard arguments and \cite[Lemma 2.11]{CE13}.  

To that end, let $\epsilon > 0$ and $L \in \bN$. Let $E_{good}(\epsilon,L)$ be the subset of $X$ such that for any $\omega \in E_{good}(\epsilon,L)$ and any $\bfv_\omega \in \bfV_\omega$, there is a subset $H(\bfv_\omega)$ of $[-1,1]$ so that \[\mu(H(\bfv_\omega)) \geq 2 - \epsilon\] 

and such that for all $s \in H(\bfv_\omega)$ \begin{align}\label{CE1}
        \lambda - \epsilon < 
            \frac{\log \sigma(g_{L} h_s, \bfv_\omega)}{L} 
            < \lambda +\epsilon 
     \end{align}

\begin{cor}\cite[Lemma 2.11]{CE13}\label{CE}
For any $\epsilon>0$ and $\delta>0$, there exists $L_0 > 0$ such that for all $L > L_0$, we have that $\nu(E_{good}(\epsilon,L)) > 1-\delta$.
\end{cor}

Similarly, for the $r_\theta$ action on a submanifold $X'$ of $X$ that supports some measure $\nu'$ with $\lambda' = \lambda(\nu')$, set $\epsilon' > 0$ and $L \in \bN$. Let $E'_{good}(\epsilon',L)$ be the subset of $X'$ such that for any $\omega \in E'_{good}(\epsilon',L)$ and any $\bfv_\omega \in \bfV_\omega$, there is a subset $H'(\bfv_\omega)$ of $[-1,1]$ so that \[\mu(H'(\bfv_\omega)) \geq 2 - \epsilon\] 

and such that for all $\theta \in H'(\bfv_\omega)$ \begin{align}\label{CE2}
        \lambda' - \epsilon' < 
            \frac{\log \sigma(g_{L} r_\theta, \bfv_\omega)}{L} 
            < \lambda' +\epsilon' 
     \end{align}

\begin{cor}\cite[Lemma 2.11]{CE13}\label{CE3}
For any $\epsilon'>0$ and $\delta'>0$, there exists $L_0 > 0$ such that for all $L > L_0$, we have that $\nu'(E'_{good}(\epsilon',L)) > 1-\delta'$.
\end{cor}

We refer to \cite[Lemma 7.5, Corollary 7.6]{AAE+17} for the adapted proofs.

\subsection{Exceptional behavior}\label{exceptional}
Since $X$ is assumed to have varying Lyapunov phenomenon, there exists an affine invariant submanifold $X'$ of $X$ that supports an ergodic $\SL_2(\bR)$-invariant measure $\nu'$, and such that $\lambda(\nu') \neq \lambda(\nu)$. Set $\lambda' = \lambda(\nu')$.  In this section, we show that in a neighborhood of the manifold $X'$, the exceptional behavior is also present in some open subset of the ambient orbit closure $X$ up to some prescribed time $L$.

For $\beta > 0$, let $X_\beta$ be the $\beta$-thick part of $X$, which is the subset of $X$ such that for all $\omega \in X$, saddle
connections have $\omega$-length at least $\beta$. For the pseudo-metric $d_{Teich}$ on $X$, define the dynamical metric

\begin{align*}
 d^n_{Teich}(x,y) \defeq \sup_{i\in[0,n]} d_{Teich}(g_n x, g_n y).
\end{align*}

Let $E'_{good}(\epsilon',L) \subset X'$ be as in Corollary \ref{CE3}. By the tubular neighborhood theorem, for any $L>0$, there exists a $\beta>0$ such that \[N_{L,\epsilon'}(X') = \left\{ \omega \in X_\beta ~:~ d^L_{Teich}(\omega,E'_{good}(\epsilon',L))< \frac{1}{L}\right\} \] is a non-empty open set that is contained in some thick part $X_\beta$ of $X$.
 
For any $\omega \in N_{L,\epsilon'}(X')$, it follows by construction that there exists $\omega' \in X'$ such that $d_{Teich}(\omega,\omega')< 1/L$ and $ d_{Teich}(g_t\omega,g_t\omega') <1/L$ for any $t\leq L$. Let $\theta_1 \in H'(\bfv_{\omega'})$ be an admissible direction as in Corollary \ref{CE3}. For all $\bfv_\omega$ in the symplectic orthogonal of $[\omega]$, we have by \ref{lipschitz} and the cocycle property that
 \begin{align*}
\log\sigma(g_t,\bfv_\omega) &\leq d_{Teich}(\omega,\omega') + \log\sigma(r_{\theta_1},\bfv_{\omega'}) + \log\sigma(g_t,r_{\theta_1}\bfv_{\omega'}) \\
&  \log\sigma(r_{\theta_2} ,g_t r_{\theta_1}\bfv_{\omega'}) + d_{Teich}(g_t\omega,g_t\omega') \\
&= d_{Teich}(\omega,\omega')  +  \log\sigma(g_t,r_{\theta_1}\bfv_{\omega'}) + d_{Teich}(g_t\omega,g_t\omega') 
 \end{align*}
In particular, we have
 \begin{align*}
\log\sigma(g_t,\bfv_\omega) &\leq d_{Teich}(\omega,\omega') + \log\sigma(g_t,r_{\theta_1}\bfv_{\omega'}) + d_{Teich}(g_t\omega,g_t\omega') \\
 &< \frac{1}{L} + (\lambda'+\epsilon')t + \frac{1}{L}  \\
 &< 2 + (\lambda'+\epsilon')t 
 \end{align*}
Similarly, we also have
 \begin{align*}
\log\sigma(g_t,\bfv_\omega) &\geq -d_{Teich}(\omega,\omega') + \log\sigma(g_t,r_{\theta_1}\bfv_{\omega'}) - d_{Teich}(g_t\omega,g_t\omega') \\
  &> -\frac{1}{L} + (\lambda'-\epsilon')t -\frac{1}{L}  \\
    &>  -2 + (\lambda'-\epsilon')t  
 \end{align*}
  for all $t\leq L$. This gives in particular that 
\begin{align}\label{123}  -2 +(\lambda'-\epsilon')t  \leq \log\sigma(g_t,\bfv_\omega) \leq 2 + (\lambda'+\epsilon')t   \end{align} for all $t\leq L$ and $\omega \in N_{L,\epsilon'}(X')$.

\subsection{Quantitative recurrence}\label{multiple}
Recall that $X = \overline{\SL_2(\bR) \omega}$. The major ingredient in our work is 

\begin{thm}\cite[Theorem 2.10]{EMM15}\label{emmhoro}
For $f \in C_c(X)$, any $\epsilon > 0$, and any interval $I \subset \bR$, there exists $T_0 > 0$ such that for all $T > T_0$, we have \[\left|\frac{1}{T} \int_0^T \frac{1}{|I|} \int_{I} f(g_t h_s \omega) ds dt - \int_{X} f d\nu \right|<\epsilon\]
\end{thm}
 
The following proposition then follows from Theorem \ref{emmhoro}, whose analogue for the $r_\theta$ action appears in the work of Chaika-Lindsey \cite[Proposition 8]{CL17}.

\begin{prop}\cite{CL17}\label{el}
Let $U_1, U_2$ be any open sets in $X$, and let $Z$ be any interval. For any $\epsilon > 0$, there exists arbitrarily large times $T>0$ such that for each $i\in\{1,2\}$, \[|\{ s \in Z ~:~ g_T h_{s} \omega \in U_i\}| \geq |Z|(\nu(U_i) - \epsilon)\]
\end{prop}

\section{Proofs of Theorems \ref{maintheta} and \ref{mainhoro}}
\subsection{Theorem \ref{mainhoro} for horocyclic arcs}
Recall we want to show that the set \[E= \left\{s \in [-1,1] ~:~ \lim_{t\to\infty} \frac{\sigma(g_{t} h_s, \bfv_\omega)}{f(t)} \text{ converges to a non-zero number }\right\}\] has zero measure with respect to the Lebesgue measure on $[-1,1]$.

We argue by contradiction. Assume that the set $E$ has positive measure. Restrict $E$ to a smaller set (which we continue to call $E$) where we have uniform convergence to the limit $C_s$ (by Egorov).  More precisely, by uniformity, we have that for any $\eta > 0$, there exists $t_\eta$ such that \[\left|\frac{\sigma(g_{t} h_s, \bfv_\omega)}{f(t)} - C_s\right| < \eta\] for all $s \in E$ for all $t > t_\eta$. This implies that for $s_1,s_2 \in E$, 
\[\frac{C_{s_1} -\eta}{C_{s_2}+\eta} < \frac{\sigma(g_{t} h_{s_1}, \bfv_\omega)}{\sigma(g_{t} h_{s_2}, \bfv_\omega)} < \frac{\eta + C_{s_1}}{C_{s_2}-\eta}\]

Observe that $C_s$ can be made continuous by applying Luzin's theorem and further restricting $E$ to a smaller set  (which we continue to call $E$), so that $C_s$ is bounded over $s \in [-1,1]$. Therefore, we have that the limits $C_s$ are bounded from above and below by $C$ and $c$, respectively. This gives us that for all $s_1,s_2 \in E(\eta)$, and all $t>t_\eta$.

\begin{equation}\label{ineq}
\frac{c -\eta}{C+\eta}  < \frac{C_{s_1} -\eta}{C_{s_2}+\eta} < \frac{\sigma(g_{t} h_{s_1}, \bfv_\omega)}{\sigma(g_{t} h_{s_2}, \bfv_\omega)} < \frac{\eta + C_{s_1}}{C_{s_2}-\eta} \leq \frac{\eta + C}{c-\eta}\end{equation}

Pick $0< \eta < c$, and let \begin{align}\label{eta}K= K(\eta) \defeq \frac{\eta + C}{c-\eta}\end{align}  so that \ref{ineq} reduces to \begin{equation}\label{ineq2}
\frac{1}{K}  <  \frac{\sigma(g_{t} h_{s_1}, \bfv_\omega)}{\sigma(g_{t} h_{s_2}, \bfv_\omega)} <  K\end{equation}
The contradiction is reached if we can show that there are arbitrary large times such that \[\frac{\sigma(g_{t} h_{s_1}, \bfv_\omega)}{\sigma(g_{t} h_{s_2}, \bfv_\omega)} > K \text { or } \frac{\sigma(g_{t} h_{s_1}, \bfv_\omega)}{\sigma(g_{t} h_{s_2}, \bfv_\omega)} < \frac{1}{K} \]  for pairs $(s_1, s_2) \in E \times E$. 

Since we assume $E$ has positive measure, there is an interval $Z = Z(\gamma)$ of some Lebesgue density point in $E$, so that \begin{equation}\label{density}\frac{|E \cap Z|}{|Z|} > 1 - \gamma\end{equation} for any $1>\gamma>0$. This is an important reduction, as it gives us access to Corollary \ref{el}, which requires an interval.

Note that although all pairs coming from $E\times E$ satisfy \ref{ineq2} for the sake of contradiction, this need not be the case for all pairs in $Z \times Z$. So we need to estimate an exceptional set of pairs, which is $ Z \times (Z - E) \cup (Z-E) \times Z$. It is immediate from \ref{density} that \begin{equation}\label{immediate}|Z - E| \leq \gamma |Z|\end{equation} From \ref{immediate}, it follows that 
\begin{equation}\label{badsetestimate}
| Z \times (Z - E) \cup (Z-E) \times Z| \leq 2\gamma |Z|^2
\end{equation}

Pick $\epsilon$ and $\epsilon'$ such that $\epsilon + \epsilon' < |\lambda - \lambda'|$, and let $\delta = \delta' = \frac{1}{2}$. Then by an application of Corollary \ref{CE} on $X$ (resp., Corollary \ref{CE3} on $X'$), there exists $L_0(\epsilon,\delta)$ (resp., $L'_0(\epsilon',\delta')$) such that the conclusion of the corollary is satisfied with $\nu(E_{good}(\epsilon, L))>1/2$ (resp., $\nu(E'_{good}(\epsilon', L))>1/2$). Pick \[L > \max\left\{L_0,L'_0, \frac{2+2\log K}{\lambda-\lambda'-\epsilon-\epsilon'},\frac{2+2\log K}{\lambda'-\lambda-\epsilon-\epsilon'}\right\}.\] 

To simplify notation, set \begin{align*}
U_1 &= E_{good}(\epsilon, L)\\
U_2 &=  N_{L,\epsilon'}(E'_{good}(\epsilon', L)).
\end{align*}

Let
\begin{align}
A_1 &=  \{s \in Z ~:~ g_t h_s \omega \in U_1\}\label{setA}\\
A_2 & = \{ s \in Z ~:~ g_t h_s \omega \in U_2\}\label{setB}
\end{align}

To show that there are pairs coming from $A_1 \times A_2$ that intersect $E \times E$, it suffices to show by \ref{badsetestimate} that $\gamma$ can be chosen so that \begin{align}
|A_1 \times A_2| > 2\gamma |Z|^2 \label{goodcond}
\end{align} Now, to ensure that \ref{goodcond} is satisfied, we want to show
\begin{align*} 
|A_1| &> \sqrt{2\gamma}|Z| + \epsilon |A_1| \implies |A_1| > \sqrt{2\gamma} \frac{|Z|}{1-\epsilon}  \\
|A_2| &> \sqrt{2\gamma}|Z| 
\end{align*}
and we note that $|A_1|$ needs to exceed $\epsilon |A_1|$ since that is the proportion of bad futures in Corollary \ref{CE} applied on $E_{good}(\epsilon, L) \subset X$. 

Choose $\gamma>0$ such that 
\begin{align}
\frac{\nu(U_1)}{2} &> \frac{\sqrt{2\gamma}}{1-\epsilon'} \\
\frac{\nu(U_2)}{2} &> \sqrt{2\gamma}
\end{align}

Let $m = \min\{\nu(U_1), \nu(U_2)\}$. By Proposition \ref{el} (with $\epsilon = m/2$), for each $i = \{1,2\}$, there exists arbitrarily large $T>t_{\eta}$ such that \begin{align}\label{mainest}
\frac{|A_i|}{|Z|} = \frac{1}{|Z|}|\{ s \in Z ~:~ g_T h_{s} \omega \in U_i\}| > \nu(U_i) - \frac{1}{2}m \geq \frac{1}{2} \nu(U_i) \end{align} 

Define 
\begin{align*}
\tau = \inf&\{ t_{\eta} \leq r\leq T : g_r h_{s_i}\omega \in  U_i \text{ for } i = \{1,2\} \\ &\text{ and for some } (s_1,s_2) \in (A_1 \times A_2) \cap (E \times E)\}
\end{align*}

Suppose $\lambda > \lambda'$. If \ref{ineq2} is not violated at $\tau$, then we have \[\frac{\sigma(g_{\tau+L} h_{s_1}, \bfv_\omega)}{\sigma(g_{\tau+L} h_{s_2}, \bfv_\omega)}\geq\frac{\exp((\lambda - \epsilon)L)\sigma(g_{\tau} h_{s_1}, \bfv_\omega)}{\exp(2+ (\lambda'+\epsilon')L)\sigma(g_{\tau} h_{s_2}, \bfv_\omega)} \geq \frac{\exp((\lambda - \epsilon)L)}{\exp(2+ (\lambda'+\epsilon')L)}\frac{1}{K}\] where we have applied the RHS of \ref{123} in the first inequality. Observe then for the chosen $L$, we have that \[\frac{\exp((\lambda - \epsilon)L)}{\exp(2+ (\lambda'+\epsilon')L)}\frac{1}{K} > K\] so that the RHS of \ref{ineq2} is violated, which is our contradiction since $(s_1,s_2)$ belongs to $E \times E$.

Now suppose instead that $\lambda < \lambda'$. If \ref{ineq2} is not violated at $\tau$, then we have \[\frac{\sigma(g_{\tau+L} h_{s_1}, \bfv_\omega)}{\sigma(g_{\tau+L} h_{s_2}, \bfv_\omega)}\leq\frac{\exp((\lambda + \epsilon)L)\sigma(g_{\tau} h_{s_1}, \bfv_\omega)}{\exp(-2+(\lambda'-\epsilon')L)\sigma(g_{\tau} h_{s_2}, \bfv_\omega)} \leq \frac{\exp((\lambda + \epsilon)L)}{\exp(-2+(\lambda'-\epsilon')L)}K\] where we have applied the LHS of \ref{123} in the first inequality. Observe then for the chosen $L$, we have \[\frac{\exp((\lambda + \epsilon)L)}{\exp(-2+(\lambda'-\epsilon')L)}K < \frac{1}{K}\] so that the LHS of \ref{ineq2} is violated, which is our contradiction since $(s_1,s_2)$ belongs to $E \times E$. \hfill\qedsymbol

\subsection{Theorem \ref{maintheta} for circle arcs}
It is straightforward to see how Theorem \ref{maintheta} can be deduced from Theorem \ref{mainhoro}. Indeed, for any $\theta \neq \pm \pi/2$, we have \[ r_\theta = \bar{h}_{\tan\theta} g_{\log\cos\theta} h_{-\tan\theta}. \]

Since $g_t\bar{h}_{\tan\theta}=\bar{h}_{e^{-2t}\tan\theta}g_{t}$, we also have  \[ g_t r_\theta = \bar{h}_{e^{-2t}\tan\theta} g_{t+\log\cos\theta} h_{-\tan\theta}. \]
So that $s$ belongs to the set in Theorem \ref{mainhoro} iff $\theta$ belongs to the set in Theorem \ref{maintheta}. \hfill\qedsymbol  

\bibliography{mybib}

\newcommand{\etalchar}[1]{$^{#1}$}
\providecommand{\bysame}{\leavevmode\hbox to3em{\hrulefill}\thinspace}
\providecommand{\MR}{\relax\ifhmode\unskip\space\fi MR }
\providecommand{\MRhref}[2]{%
  \href{http://www.ams.org/mathscinet-getitem?mr=#1}{#2}
}
\providecommand{\href}[2]{#2}
\begin{thebibliography}{ASAE{\etalchar{+}}17}

\bibitem[AN20]{AN19}
David Aulicino and Chaya Norton, \emph{Shimura-{T}eichm\"{u}ller curves in
  genus 5}, J. Mod. Dyn. \textbf{16} (2020), 255--288. \MR{4160181}

\bibitem[AS21]{A19}
Hamid Al-Saqban, \emph{{A Central Limit Theorem for the Kontsevich-Zorich
  cocycle}}, preprint (2021).

\bibitem[ASAE{\etalchar{+}}17]{AAE+17}
Hamid Al-Saqban, Paul Apisa, Alena Erchenko, Osama Khalil, Shahriar Mirzadeh,
  and Caglar Uyanik, \emph{{Exceptional directions for the Teichm\"uller
  geodesic flow and Hausdorff dimension}}, arXiv preprint arXiv:1711.10542 (to
  appear in Journal of Eur. Math. Soc.). (2017).

\bibitem[Aul16]{Aul16}
David Aulicino, \emph{{Affine invariant submanifolds with completely degenerate
  Kontsevich--Zorich spectrum}}, Ergodic Theory and Dynamical Systems (2016),
  1--24.

\bibitem[BSW16]{BSW16}
Matt Bainbridge, John Smillie, and Barak Weiss, \emph{Horocycle dynamics: new
  invariants and eigenform loci in the stratum {H}{(1,1)}}, arXiv preprint
  arXiv:1603.00808 (to appear in Memoirs of the AMS) (2016).

\bibitem[CE15]{CE13}
Jon Chaika and Alex Eskin, \emph{Every flat surface is {B}irkhoff and
  {O}seledets generic in almost every direction}, J. Mod. Dyn. \textbf{9}
  (2015), 1--23. \MR{3395258}

\bibitem[CKS21]{CKS21}
Jon Chaika, Osama Khalil, and John Smillie, \emph{{On the space of ergodic
  measures for the horocycle flow on strata of Abelian differentials}},
  preprint (2021).

\bibitem[CL17]{CL17}
Jon Chaika and Kathryn Lindsey, \emph{{Horocycle flow orbits and lattice
  surface characterizations}}, Ergodic Theory and Dynamical Systems (2017),
  1--21.

\bibitem[CM13]{CM13}
Dawei Chen and Martin M{\"o}ller, \emph{Nonvarying sums of {L}yapunov exponents
  of {A}belian differentials in low genus}, Geometry \& Topology \textbf{16}
  (2013), no.~4, 2427--2479.

\bibitem[DFV17]{DFV17}
Dmitry Dolgopyat, Bassam Fayad, and Ilya Vinogradov, \emph{Central limit
  theorems for simultaneous diophantine approximations}, Journal de
  l’{\'E}cole polytechnique-Math{\'e}matiques \textbf{4} (2017), 1--35.

\bibitem[EKZ14]{EKZ14}
Alex Eskin, Maxim Kontsevich, and Anton Zorich, \emph{{Sum of Lyapunov
  exponents of the Hodge bundle with respect to the Teichm{\"u}ller geodesic
  flow}}, Publications math{\'e}matiques de l'IH{\'E}S \textbf{120} (2014),
  no.~1, 207--333.

\bibitem[EM18]{EM18}
Alex Eskin and Maryam Mirzakhani, \emph{{Invariant and stationary measures for
  the SL(2,R) action on moduli space}}, Publications math{\'e}matiques de
  l'IH{\'E}S \textbf{127} (2018), no.~1, 95--324.

\bibitem[EMM15]{EMM15}
Alex Eskin, Maryam Mirzakhani, and Amir Mohammadi, \emph{{Isolation,
  equidistribution, and orbit closures for the SL(2,R) action on moduli
  space}}, Annals of Mathematics (2015), 673--721.

\bibitem[Fil16]{Fil06}
Simion Filip, \emph{{Semisimplicity and rigidity of the Kontsevich-Zorich
  cocycle}}, Inventiones mathematicae \textbf{205} (2016), no.~3, 617--670.

\bibitem[FM08]{FM08}
Giovanni Forni and Carlos Matheus, \emph{{An example of a Teichmuller disk in
  genus 4 with degenerate Kontsevich-Zorich spectrum}}, arXiv preprint
  arXiv:0810.0023 (2008).

\bibitem[FM14]{FM13}
\bysame, \emph{Introduction to {T}eichm\"{u}ller theory and its applications to
  dynamics of interval exchange transformations, flows on surfaces and
  billiards}, J. Mod. Dyn. \textbf{8} (2014), no.~3-4, 271--436. \MR{3345837}

\bibitem[FMZ11]{FMZ11}
Giovanni Forni, Carlos Matheus, and Anton Zorich, \emph{Square-tiled cyclic
  covers}, J. Mod. Dyn. \textbf{5} (2011), no.~2, 285--318. \MR{2820563}

\bibitem[FMZ12]{FMZ12}
\bysame, \emph{{Lyapunov spectrum of invariant subbundles of the Hodge
  bundle}}, Ergodic Theory and Dynamical Systems \textbf{34} (2012), no.~2,
  353--408.

\bibitem[For02]{For02}
Giovanni Forni, \emph{{Deviation of ergodic averages for area-preserving flows
  on surfaces of higher genus}}, Annals of Mathematics \textbf{155} (2002),
  no.~1, 1--103.

\bibitem[For06]{For06}
\bysame, \emph{{On the Lyapunov exponents of the Kontsevich--Zorich cocycle}},
  Handbook of dynamical systems \textbf{1} (2006), 549--580.

\bibitem[M{\"o}l11]{Mol05}
Martin M{\"o}ller, \emph{Shimura and {T}eichm\"{u}ller curves}, J. Mod. Dyn.
  \textbf{5} (2011), no.~1, 1--32. \MR{2787595}

\end{thebibliography}
\bibliographystyle{amsalpha}
\end{document}